\documentclass[11pt]{amsart}
\usepackage{amssymb,amscd}
\usepackage[arrow,matrix,curve]{xy}
\newcommand\RE{\mathbb{R}}
\newcommand\ZA{\mathbb{Z}}

\newcommand\TT{\mathbb{T}}

\newcommand\Tr{\operatorname{Tr}}
\newcommand\Det{\operatorname{Det}}

\newcommand\tr{\operatorname{tr}}

\newcommand\im{\operatorname{im}}
\newcommand\bra{\langle}
\newcommand\ket{\rangle}
\newcommand{\medwedge}{\wedge}

\newcommand\del{\delta}
\newcommand\lam{\lambda}
\newcommand\Gam{\mathit{\Gamma}}
\newcommand\Lam{\mathit{\Lambda}}
\newcommand\Om{\mathit{\Omega}}

\theoremstyle{plain}
\newtheorem{theorem}{Theorem}[section]
\newtheorem{lemma}[theorem]{Lemma}

\theoremstyle{definition}

\newtheorem{definition*}{Definition}

\theoremstyle{remark}

\newtheorem{remarks*}{Remarks}

\begin{document}

\begin{flushright}
{\tt arXiv:0912.2184v2[math.DG]}\\
February, 2010 (revised)
\end{flushright}

\vspace{1cm}

\title[Twisted Analytic Torsion]
{Twisted Analytic Torsion}

\author{Varghese Mathai}
\address{Department of Mathematics, University of Adelaide,
Adelaide 5005, Australia}
\email{mathai.varghese@adelaide.edu.au}

\author{Siye Wu}
\address{Department of Mathematics, University of Colorado,
Boulder, Colorado 80309-0395, USA and
Department of Mathematics, University of Hong Kong, Pokfulam Road,
Hong Kong, China}
\email{swu@math.colorado.edu\quad swu@maths.hku.hk}

\begin{abstract}
\vspace{5ex}

We review the Reidemeister torsion, Ray-Singer's analytic torsion and the 
Cheeger-M\"uller theorem.
We describe the analytic torsion of the de Rham complex twisted by a flux 
form introduced by the current authors and recall its properties.
We define a new twisted analytic torsion for the complex of invariant
differential forms on the total space of a principal circle bundle twisted
by an invariant flux form. 
We show that when the dimension is even, such a torsion is invariant
under certain deformation of the metric and the flux form.
Under $T$-duality which exchanges the topology of the bundle and the
flux form and the radius of the circular fiber with its inverse, 
the twisted torsions of invariant forms are inverse to each other for
any dimension.\\

\end{abstract}

\maketitle

\vspace{-3ex}

\begin{center}
{\small Dedicated to Professor Yang Lo on the occasion of his 70th birthday}
\end{center}

\vspace{1cm}

\section*{Introduction}

Reidemeister torsion (or $R$-torsion) was introduced by Reidemeister \cite{R}
for $3$-manifolds.
It was generalized to higher odd dimensions by Franz \cite{F} and de Rham
\cite{dR36}.
As the first invariant that could distinguish spaces which are homotopic but
not homeomorphic, it can be used to classify lens spaces \cite{F,dR64,Mi}.
Analytic torsion (or Ray-Singer torsion) is a smooth invariant of compact
manifolds of odd dimensions defined by Ray and Singer \cite{RS,RS2} as an
analytic counterpart of the Reidemeister torsion.
They conjectured that the two torsions are equal for compact manifolds.
This Ray-Singer conjecture was proved independently by Cheeger \cite{C}
and M\"uller \cite{M78}.
Another proof of the Cheeger-M\"uller theorem with an extension was given
by Bismut and Zhang \cite{BZ} using Witten's deformation.
In \cite{MW}, the current authors introduced an analytic torsion for the
de Rham complex twisted by a flux form.
It shares many properties with the Ray-Singer torsion but has several novel
features. 
In this paper, we review these developments and introduce a new
twisted analytic torsion for the complex of invariant differential forms on
a principal circle bundle.

The paper is organized as follows.
In Section~1, we review Reidemeister's combinatorial torsion, Ray and
Singer's analytic torsion and the Cheeger-M\"uller theorem.
In Section~2, we describe the analytic torsion of the de Rham complex
twisted by a flux form, its invariance under the deformation of the
metric and flux form when the dimension is odd, the relation with
generalized geometry when the flux is a $3$-form and the behavior
under $T$-duality for $3$-manifolds.
In Section~3, we introduce a new twisted analytic torsion for the
complex of invariant differential forms on the total space of a
principal circle bundle twisted by an invariant flux form. 
We show that when the dimension is even, such a torsion is invariant
under certain deformation of the metric and the flux form.
Under $T$-duality, which exchanges the topology of the bundle with
the flux form and the radius of circular fibers with its inverse, 
the twisted torsions are inverse to each other for any dimensions.

\section{Reidemeister and Ray-Singer torsions}

In algebraic topology, several groups can be assigned to a topological space
such as its fundamental group, homology and cohomology groups.
For example, the dimension of $p$-th cohomology group $H^p(X,\RE)$ is,
roughly speaking, the number of $p$-dimensional ``holes'' in the space $X$.
These groups are invariants in the sense that if two topological spaces
are the same, or more precisely homeomorphic, then the corresponding
groups are isomorphic.
So if we find two spaces $X$ and $Y$ with $H^p(X,\RE)\not\cong H^p(Y,\RE)$ for
some integer $p$, then we can conclude that $X$ is not homeomorphic to $Y$.
In this way, we can distinguish different topological spaces.

However, these groups are invariant not just under homeomorphisms,
but under continuous deformations, or homotopy.
Two spaces can be homotopic but not homeomorphic to each other.
As none of the above groups can distinguish two homotopic yet different
topological spaces, we need additional invariants to achieve this.
Reidemeister torsion is such an invariant that is invariant under
homeomorphisms but not under homotopy.
It is a secondary invariant in the sense that it is an element in a certain
space constructed from the above groups (which are regarded as the primary
invariants).

To recall its definition, we consider a topological space of dimension $n$
with a finite triangulation $K$.
The simplicial structure gives rise to a cochain complex $(C^\bullet(K),\del)$;
we will take the real coefficients unless otherwise indicated.
Roughly, the Reidemeister torsion is the alternating product of the
determinants of the coboundary operators $\delta$ on various $C^p(K)$
($p=0,1,\dots,n$).
More precisely, we take the adjoint $\delta^\dagger$ of
$\delta\colon C^p(K)\to C^{p+1}(K)$ with respect to the inner product under
which the $p$-simplices form an orthonormal basis.
We define the Laplacians $\Delta_p=\delta^\dagger\delta+\delta\delta^\dagger$
on $C^p(K)$.
Just like in Hodge theory, we have $H^p(K)\cong\ker\Delta_p$.
If $H^p(K)=0$ were true for any $p$, then all $\Delta_p$ would be invertible,
and we could define the Reidemeister torsion as the number \cite{R,F,dR36}
\[   \tau(K)=\prod_{p=0}^n(\det\Delta_p)^{(-1)^{p+1}p/2}.     \]
To take into account of the non-trivial cohomology groups, we choose a unit
volume element $\eta_p$ of $\ker\Delta_p\subset C^p(K)$ for each $p$.
Thus, we have $\eta_p\in\det\ker\Delta_p\cong\det H^p(K)$.
We define the Reidemeister torsion 
\[  \tau(K)=\prod_{p=0}^n(\det{}'\Delta_p)^{(-1)^{p+1}p/2}
        \bigotimes_{p=0}^n\eta_p^{(-1)^p}                    \]
as an element of the line 
$\det H^\bullet(K)=\bigotimes_{p=0}^n\det H^p(K)^{(-1)^p}$
constructed from the cohomology groups.
Here, $\det'$ means taking the determinant of an operator in the subspace
orthogonal to its kernel. 

We make two remarks here.
First, in the above construction, the torsion is only defined up to a sign as
the unit volume elements depend on the orientation on the cohomology groups.
A more intrinsic way is to define torsion as a norm on the determinant line
bundle \cite{BZ}.
Second, the cochains on $K$ can take values in a local system that comes
from an orthogonal or unitary representation of the fundamental group of
$K$ or equivalently, a flat vector bundle over the underlying
topological space.
The torsion defined would then depend on this data.

It can be shown that $\tau(K)$ is invariant under the subdivision of $K$
and is invariant under homeomorphism.
It also satisfy some functorial properties.
Happily, it is not invariant under homotopy.
A celebrated example is the lens space, which is the quotient of
the $3$-sphere by a finite Abelian group.
In fact, the lens spaces are classified by their cohomology groups
and the Reidemeister torsion \cite{F,dR64,Mi}.
More recently, equivariant torsion has been used to classify isometries and
quotients of certain symmetric spaces up to diffeomorphisms \cite{Ko}.

An alternative way to define the cohomology groups is by differential forms.
If $X$ is a smooth manifold, let $d\colon\Om^p(X)\to\Om^{p+1}(X)$
be the exterior differentiation of the space of $p$-forms, with $d^2=0$.
The $p$-th de Rham cohomology is
\[ H_{\mathrm{dR}}^p(X)=\frac{\ker(d\colon\Om^p(X)\to\Om^{p+1}(X))}
                         {\im(d\colon\Om^{p-1}(X)\to\Om^p(X))}.    \]
The de Rham theorem states that there is a natural isomorphism 
$H_{\mathrm{dR}}^p(X)\cong H^p(K)$, where $K$ is a triangulation of $X$.

A natural question is how to represent the Reidemeister torsion analytically.
Suppose $X$ is compact, orientable and is equipped with a Riemannian metric.
Then we define the Laplacians $\Delta_p=d^\dagger d+dd^\dagger$,
where the adjoint is with respect to the inner product on $\Om^p(X)$
given by the metric.
If $X$ is compact, then $H^p_{\mathrm{dR}}(X)\cong\ker\Delta_p$ and is
finite dimensional by Hodge theory.
Taking the unit volume form $\eta_p$ of $\ker\Delta_p\subset\Om^p(X)$, 
the Ray-Singer analytic torsion is
\[ \tau_{\mathrm RS}(X)=\prod_{p=0}^n(\Det'\Delta_p)^{(-1)^{p+1}p/2}
        \bigotimes_{p=0}^n\eta_p^{(-1)^p}                      \]
as an element, defined up to a sign, of the line
$\det H^\bullet_{\mathrm{dR}}(X)$ \cite{RS}.
The above construction generalizes easily if the differential forms are
valued in a flat Hermitian vector bundle.
By replacing $d$ in the Laplacian with the flat connection, we obtain 
Ray-Singer torsion that depends on a flat vector bundle.

The determinants of the Laplacians, and of many other operators on infinite
dimensional spaces, are defined by regularization using zeta-functions.
Let $A$ be a self-adjoint semi-positive operator acting on a Hilbert space. 
Suppose the positive eigenvalues are listed as $0<\lam_1\le\lam_2\le\cdots$,
taking into account of multiplicities.
Since the eigenvalues $\lam_i$ typically go to infinity as $i\to\infty$, it
does not make sense to consider their product in the usual sense.
We define the zeta-function of $A$ as
\[   \zeta_A(s)=\sum_{i=1}^\infty\lam_i^{-s}=\tr' A^{-s},   \]
which is a function of one complex variable $s$.
Here $\tr'$ means taking the trace in the subspace orthogonal to $\ker A$. 
If $A$ is elliptic, then by the heat kernel method, $\zeta_A(s)$ is analytic
in $s$ when $\Re s$ is sufficiently large and it can be extended 
meromorphically to the complex plane so that it is analytic at $s=0$.
We define
\[    \Det'A=e^{-\zeta'_A(0)}.    \]
When $A$ acts on a finite dimensional space, the above definition reduces
to the finite product of the positive eigenvalues of $A$.

The Ray-Singer torsion is invariant under the deformation of the Riemannian
metric on $X$ when $\dim X$ is odd and therefore is a smooth invariant of 
odd-dimensional compact manifolds.
Furthermore, it satisfies the same functorial properties as the Reidemeister
torsion.
Ray and Singer \cite{RS,RS2} therefore conjectured that their torsion is 
equal to the Reidemeister torsion for compact odd-dimensional manifolds.
The case for lens spaces was verified in \cite{Ray}.
The conjecture was proved independently by Cheeger \cite{C} and
M\"uller \cite{M78}, when the local coefficients are given by an orthogonal or
unitary representation of the fundamental group.
The result is also true when the representation is unimodular \cite{M93}.
In \cite{BZ}, Bismut and Zhang extended the Cheeger-M\"uller theorem
using the deformation of Witten \cite{Wi}.
They established the variation of the torsion under the deformations of
the metric on the manifold of arbitrary dimensions and the Hermitian form 
on the vector bundle with an arbitrary flat connection.

\section{Analytic torsion of twisted complexes}

Let $X$ be a smooth manifold and $H$, a closed $3$-form $H$ called flux form,
which has its origin in supergravity and string theory \cite{RW,AS}.
We have an operator $d_H=d+H\wedge\cdot\,$ acting on $\Om^\bullet(X)$, which
squares to zero.
So it can be used to define a twisted de Rham cohomology $H^\bullet(X,H)$.
The twisted de Rham complex is only $\ZA_2$-graded and so is the twisted
de Rham cohomology.
Let $\Om^{\bar0}(X)$, $\Om^{\bar1}(X)$ denote the space of differential 
forms on $X$ of even, odd degrees, respectively.
Then the operator $d_H$ acts between these two spaces.

In fact, the above flux form can be a closed form of any odd degree
with no $1$-form component.
In \cite{MW}, we introduced the analytic torsion of the operator $d_H$
when $X$ is compact, which we now assume.
The main difficulty was that the twisted de Rham complex does not have
a $\ZA$-grading.
Formally, given a Riemannian metric on $X$, the twisted analytic torsion is
\[ \tau(X,H)=\Det'_{\Om^{\bar0}(X)}(d_H^\dagger d_H)^{1/2}
             \Det'_{\Om^{\bar1}(X)}(d_H^\dagger d_H)^{-1/2}
             \,\eta_{\bar0}\otimes\eta_{\bar1}^{-1},                 \]
where $\eta_{\bar k}$ is a unit volume element of $H^{\bar k}(X,H)$ ($k=0,1$).
However, unlike the Laplacian $\Delta_H=d_H^\dagger d_H+d_Hd_H^\dagger$, the
operator $d_H^\dagger d_H$ on $\Om^{\bar k}(X)$ is not elliptic and the usual
heat kernel methods seem inadequate.
Instead, we use some properties of pseudo-differential operators.
Let $P_{\bar k}$ be the projection in $\Om^{\bar k}(X)$ onto the 
image of $d_H^\dagger$.
Since $P_{\bar k}$ is a pseudo-differential operator, we have \cite{Se}
\[  \tr(P_{\bar k}\Delta_H^{-s})=\frac{c_{-1}}{s}+c_0+o(s),   \]
upon meromorphically continuing the left-hand side to $s=0$. 
The coefficient $c_{-1}$ is essentially the non-commutative residue
\cite{Wo,Gui} of $P_{\bar k}$.
It turns out that since $P_{\bar k}$ is also a projection, 
$c_{-1}=0$ \cite{Wo}.
Consequently, the zeta-function of the restriction of $d_H^\dagger d_H$ to the
image of $d_H^\dagger$ is regular at $s=0$ and its determinant is defined.
It is also possible to include a flat Hermitian vector bundle in the
definition of twisted analytic torsion.

Just as the usual de Rham cohomology groups, the twisted counterpart
$H^\bullet(X,H)$ is also invariant under homotopies \cite{MW}.
The twisted torsion $\tau(X,H)$ satisfies a similar set of functorial
properties as in \cite{RS}.
Moreover, when $X$ is compact and $\dim X$ is odd, $\tau(X,H)$ is invariant
under the deformation of the Riemannian metric on $X$, the inner product
or Hermitian structure on the flat bundle and the flux form within its 
cohomology class \cite{MW}. 
Therefore the analytic torsion $\tau(X,H)$ is a secondary invariant in
the same sense but in the twisted setting. 

When $H$ is a $3$-form on $X$, the deformation of the Riemannian metric $g$
and that of the flux form $H$ within its cohomology class can be interpreted
as a deformation of generalized metrics on $X$ \cite{MW}.
Recall that a generalized metric on $X$ is a splitting 
$TX\oplus T^*X=T_+X\oplus T_-X$ such that the bilinear form
\[   \bra x+\alpha,y+\beta\ket=(\alpha(y)+\beta(x))/2,      \]
where $x,y$ are vectors fields and $\alpha,\beta$ are $1$-forms on $X$,
is positive definite (negative definite, respectively) on $T_\pm X$ and
such that $\bra T_+X,T_-X\ket=0$ \cite{Gua2}.
The subspace $T_+X$ is the graph of $g+B$, where $g$ is a usual Riemannian
metric and $B$ is a $2$-form. 
Given a generalized metric, there is an inner product on $\Om^\bullet(X)$
called the Born-Infeld metric \cite{Gua2}. 
It can be shown that the effect of deformation $H\mapsto H-dB$ on torsion
is equivalent to taking the adjoint of the operator $d_H$ with respect to
the Born-Infeld metric \cite{MW}.
This amounts to deforming a usual metric $g$ to a generalized one.
Thus analytic torsion should be defined for generalized metrics so that the
deformations of $g$ and of $H$ are unified.

We consider a special case when $H$ is a top form on an odd-dimensional
orientable compact manifold $X$.
The cohomology class of $H$ is $[H]\in H^{\mathrm{top}}(X,\RE)\cong\RE$. 
Then the twisted analytic torsion is \cite{MW}
\[  \tau(X,H)=|[H]|\,\tau_{\mathrm{RS}}(X)\,
    \eta_{\bar0}\otimes\eta_{\bar1}^{-1},  \] 
where $\tau_{\mathrm{RS}}(X)$ is the Ray-Singer torsion.
This provides examples of twisted torsion for $3$-manifolds, when $H$
has to be a $3$-form. 
The calculation of $\tau(X,H)$ is based on the work of Kontsevich and Vishik
\cite{KV} on factorization of determinants in odd dimensions, although we hope
that there are simpler methods\footnote{For example, affirmative answers to
the following questions on functions of one complex variable would provide a
simpler calculation.
Suppose a sequence of positive real numbers $0<\lam_1\le\lam_2\le\cdots$
goes to infinity fast enough so that the zeta-function
$\zeta(s)=\sum_{i=1}^\infty\lam_i^{-s}$ is absolutely convergent and hence
analytic when $\Re s$ is sufficiently large. Suppose $\zeta(s)$ can be
meromorphically continued so that it is regular at $s=0$.
We partition the set of positive integers to a disjoint union of finite
sets $I_k$ ($k=1,2,\dots$).
Let $\Lam_k=\prod_{i\in I_k}\lam_i$.
Then $Z(s)=\sum_{k=1}^\infty\Lam_k^{-s}$ is absolutely convergent and hence
analytic when $\Re s$ is sufficiently large.
Can $Z(s)$ be meromorphically continued so that it is regular at $s=0$?
If so, is it true that $Z'(0)=\zeta'(0)$?}
which can be useful in more general cases as well.
The factor $|[H]|$ is also the torsion of the spectral sequence in \cite{RW}.
It is not clear whether there is a simple relation between the twisted
torsion and the classical Ray-Singer torsion in the general situation.

Whereas Reidemeister's combinatorial torsion precedes Ray-Singer's analytic
torsion, the simplicial counterpart of the twisted analytic torsion is
still missing in the general case.
This is because the cup product in the simplicial cochain complex is
associative but not graded commutative.
Let $h$ be a simplicial cocycle that represents the same cohomology class
(under the de Rham isomorphism) as $H$.
Then $(\delta+h\cup\,\cdot\,)^2=(h\cup h)\cup\,\cdot\,$ may not be zero.
However, the situation simplifies if the degree of $H$ or $h$ is greater
than $\dim X/2$, in which case $h\cup h=0$ for dimension reason. 
In particular, this condition is satisfied if $H$ is a top-degree form.
In fact, a twisted version of the Cheeger-M\"uller theorem holds in this case
\cite{MW}.

We consider the behavior of the twisted torsion under $T$-duality.
Let $\TT$ be the circle group and $\pi\colon X\to M$, a principal $\TT$-bundle
with Euler class $e(X)\in H^2(M,\ZA)$.
Suppose $H$ is a closed $3$-form on $X$ with integral periods.
By the Gysin sequence, there is a dual circle fibration 
$\hat\pi\colon \hat X\to M$, whose fiber is the dual circle group $\hat\TT$, 
with a flux $3$-form $\hat H$ on $\hat X$
such that \cite{BEM} 
\[    \pi_*[H]=e(\hat X),\quad \hat\pi_*[\hat H]=e(X)   \]
(modulo torsion elements). 
Thus $T$-duality for circle bundles is the exchange of background $H$-flux
on the one side and the Chern class on the other.
We have the following duality on twisted cohomology groups \cite{BEM}: 
\[   H^{\bar0}(\hat X,\hat H)\cong H^{\bar1}(X,H),\quad
     H^{\bar1}(\hat X,\hat H)\cong H^{\bar0}(X,H).       \]    
Consequently, $\det H^\bullet(\hat X,\hat H)\cong(\det H^\bullet(X,H))^{-1}$.
When $\dim X=3$, since $H$ is a top-degree form, we get \cite{MW}
\[  \tau(\hat X,\hat H)=\tau(X,H)^{-1}   \]
under the above identification.
In general, it is not known whether there is such a concise relation.
Instead, we will explore an invariant version of the twisted torsion
and its behavior under $T$-duality in the next section.

Our method applies to other $\ZA_2$-graded complexes \cite{MW2}.
For example, suppose $X$ is a complex manifold and $H\in\Om^{0,\bar1}(X)$.
Let $\bar\partial_H=\bar\partial+H\wedge\cdot\,$.
If $\bar\partial H=0$, then $\bar\partial_H^2=0$ just as in the de Rham case.
Using the same argument, we can define an analytic torsion of the twisted
Dolbeault complex \cite{MW2}.
Alternatively, we can take $H\in\Om^{1,2}(X)$ with $\bar\partial H=0$.
Let $\Om^{-p,q}(X)=\Gam(\medwedge^pT^{1,0}(X)\otimes\medwedge^q(T^{0,1}X)^*)$.
Then $H\wedge\cdot\,\colon\Om^{-p,q}(X)\to\Om^{-(p-1),q+2}(X)$ and again,
$\bar\partial_H^2=0$.
The cohomology of $\bar\partial_H$ contains information of the deformation of
twisted generalized complex structures \cite{Li}.
A torsion as a secondary invariant in this case can also be defined.

\section{$T$-duality for circle bundles and analytic torsion}

\subsection{Analytic torsion for the complex of invariant forms}
Consider a principal $\TT$-bundle $\pi\colon X\to M$.
Suppose $X$ is compact and $H$ is a $\TT$-invariant closed $3$-form on $X$.
We consider the $\ZA_2$-graded complex $(\Om^\bullet(X)^\TT,d_H)$ of
$\TT$-invariant differential forms on $X$.
At first sight, it seems difficult to define torsion, as the asymptotic
expansion of the heat kernel with a group action \cite{BrHe}, and hence the
pole structure of the corresponding zeta-function, is rather complicated.
However, given a connection on $X$, the space $\Om^\bullet(X)^\TT$ is
isomorphic to $\Om^\bullet(M)\oplus\Om^\bullet(M)$.
Through this isomorphism, operators on $\Om^\bullet(X)^\TT$ acts on sections 
of bundles over $M$.
The torsion is then defined by the regularized determinants of elliptic
operators on $M$.
 
Suppose the metric $g_X$ on $X$ is $\TT$-invariant and such that the length
of every circular fiber is equal to some constant $r>0$.
Then
\[ g_X=\pi^*g_M+r^2 A\odot A, \]
where $g_M$ is a metric on $M$, $A\in\Om^1(X)$ is a connection 
$1$-form with the normalization $\int_\TT A=1$ and $\odot$ stands
for the symmetric tensor product.
Since any $\omega\in\Om^{\bar k}(X)^\TT$ can be written uniquely as
$\omega=r^{1/2}\pi^*\omega_1+r^{-1/2}A\wedge\pi^*\omega_2$ with
$\omega_1\in\Om^{\bar k}(M)$ and $\omega_2\in\Om^{\overline{k-1}}(M)$,
there is an isomorphism
$$
\phi\colon\Om^{\bar k}(X)^\TT\to\Om^{\bar k}(M)\oplus\Om^{\overline{k-1}}(M)
$$
for $k=0,1$ defined by $\phi(\omega)=(\omega_1,\omega_2)$.

\begin{lemma}\label{isom}
The above isomorphism $\phi$ is an isometry under the inner products on
$\Om^\bullet(X)^\TT$ and $\Om^\bullet(M)$ defined by $g_X$ and $g_M$,
respectively.
\end{lemma}

\begin{proof}
We have
$*_X\,\omega=r^{-1/2}A\wedge\pi^*(*_M\,\omega_1)+r^{1/2}\pi^*(*_M\,\omega_2)$
and thus
$$
\int_X\omega\wedge*_X\,\omega=
\int_M\omega_1\wedge*_M\,\omega_1+\int_M\omega_2\wedge*_M\,\omega_2.
$$
The result follows.
\end{proof}

Let $d_{\bar k}$ be the restriction of the operator $d_H$ to
$\Om^{\bar k}(X)^\TT$ and let $\tilde d_{\bar k}=\phi d_{\bar k}\phi^{-1}$.
If we write $H=\pi^*H_3+A\wedge\pi^* H_2$ with $H_3\in\Om^3(M)$ and
$H_2\in\Om^2(M)$ and denote by $F\in\Om^2(M)$ the curvature $2$-form of $A$
(that is, $\pi^*F=dA$), then
$$
\tilde d_{\bar k}={d_{H_3}\quad\;r^{-1}F\choose rH_2\quad-d_{H_3}}
$$ 
on $\Om^{\bar k}(M)\oplus\Om^{\overline{k-1}}(M)$.
Since $\phi$ is an isometry, we have 
$\tilde d_{\bar k}^\dagger=\phi d_{\bar k}^\dagger\phi^{-1}$ and
$\tilde\Delta_{\bar k}:=\tilde d_{\bar k}^\dagger\tilde d_{\bar k}+
\tilde d_{\overline{k-1}}\tilde d_{\overline{k-1}}^\dagger
=\phi\Delta_{\bar k}\phi^{-1}$.

Clearly, $\tilde\Delta_{\bar k}$ is a second order elliptic operator on $M$
whose leading symbol is the same as the (untwisted) Laplacian.
The projection $\tilde P_{\bar k}=\phi P_{\bar k}\phi^{-1}$ onto the
closure of $\im(\tilde d_{\bar k}^\dagger)$ is a pseudo-differential
operator on $M$ of degree zero.
By the same argument as in Theorem~2.1 of \cite{MW}, the zeta-function
$$
\zeta_\TT(s,d_{\bar k}^\dagger d_{\bar k}):=\Tr_{\Omega^{\bar k}(X)^\TT}
P_{\bar k}\Delta_{\bar k}^{-s}=\Tr\tilde P_{\bar k}\tilde\Delta_{\bar k}^{-s} 
$$
is holomorphic at $s=0$ and hence we can define the determinant
$$
\Det'_\TT d_{\bar k}^\dagger d_{\bar k}
:=\Det'\tilde d_{\bar k}^\dagger\tilde d_{\bar k}
=e^{-\zeta'_\TT(0, d_{\bar k}^\dagger d_{\bar k})}.
$$
The analytic torsion for the $\TT$-invariant part of the twisted de Rham
complex is then defined, up to a sign, as
$$
\tau_\TT(X,H,r):=(\Det'_\TT d_{\bar 0}^\dagger d_{\bar 0})^{1/2}
(\Det'_\TT d_{\bar 1}^\dagger d_{\bar 1})^{-1/2}
\eta_{\bar0}\otimes\eta_{\bar1}^{-1}\in\det H^\bullet(X,H),
$$
where the unit volume elements $\eta_{\bar k}$ ($k=0,1$) are as before,
since $\ker(\Delta_{\bar k})\cong H^{\bar k}(X,H)$ are invariant under $\TT$.
In fact, if $\tilde\eta_{\bar k}$ is the unit volume element of 
$\ker(\tilde\Delta_{\bar k})$, then since $\phi$ is an isometry,
$\tilde\eta_{\bar k}=(\det\phi)(\eta_{\bar k})$.
 
\begin{theorem} 
If $X$ is compact and $\dim X$ is even, then $\tau_\TT(X,H,r)$ is invariant 
under the deformations of $g_X$ in the class of $\TT$-invariant metrics
such that the length of every fiber is $r$ and the deformations of $H$ in the 
space of $\TT$-invariant $3$-forms representing the same cohomology class.
\end{theorem}

\begin{proof}
If the metric $g_M$ is deformed along a path parametrized by $u\in\RE$, then
$$
\frac{\partial\tilde d_{\bar k}}{\partial u}=0,\qquad
\frac{\partial\tilde d_{\bar k}^\dagger}{\partial u}
=-[\tilde\alpha,\tilde d_{\bar k}^\dagger],
$$
where $\tilde\alpha:=\frac{\partial(*_M)}{\partial u}$.
If the $3$-form $H$ is deformed in its cohomology class with a parameter
$v\in\RE$, let $B\in\Om^2(X)$ be given by $\frac{\partial H}{\partial v}=-dB$.
Then
$$
\frac{\partial\tilde d_{\bar k}}{\partial v}=[\tilde\beta,\tilde d_{\bar k}],
\qquad\frac{\partial\tilde d_{\bar k}^\dagger}{\partial v}=-
[\tilde\beta^\dagger,\tilde d_{\bar k}^\dagger],
$$
where $\tilde\beta=\phi\beta\phi^{-1}$ and $\beta=B\wedge\cdot\,$.
Since $\dim M$ is odd, the method of \S 3.1 and \S 3.2 of \cite{MW} shows that 
$$
(\Det'_\TT\tilde d_{\bar 0}^\dagger\tilde d_{\bar 0})^{1/2}
(\Det'_\TT\tilde d_{\bar 1}^\dagger\tilde d_{\bar 1})^{-1/2}
\tilde\eta_{\bar0}\otimes\tilde\eta_{\bar1}^{-1}
$$
is invariant under these deformations.
In both cases, the isomorphism $\phi$ remains fixed.
So $\tau_\TT(X,H,r)$ is also invariant.

Finally, if the connection $A$ is deformed with a parameter $w\in\RE$, 
let $C=-\frac{\partial A}{\partial w}$.
Then
$$
\frac{\partial\tilde d_{\bar k}}{\partial w}=[\tilde\gamma,\tilde d_{\bar k}],
\qquad\frac{\partial\tilde d_{\bar k}^\dagger}{\partial w}=
-[\tilde\gamma^\dagger,\tilde d_{\bar k}^\dagger],
$$
where $\tilde\gamma={0\quad\!\!r^{-1}C\choose 0\quad\;0\quad}$.
Since $\dim M$ is odd, following the proof of Lemma~3.5 of \cite{MW}, we get
$$
\frac{\partial}{\partial w}
\log[\Det'\tilde d_{\bar 0}^\dagger\tilde d_{\bar 0}\,
(\Det'\tilde d_{\bar 1}^\dagger\tilde d_{\bar 1})^{-1}]
=2\sum_{k=0,1}(-1)^k\Tr(\tilde\gamma\tilde Q_{\bar k}),
$$
where $\tilde Q_{\bar k}=\phi Q_{\bar k}\phi^{-1}$ is the projection onto
$\ker(\tilde\Delta_{\bar k})$.
On the other hand, the isomorphism $\phi$ varies as $w$.
If $\phi^{(0)}$ is the isomorphism at $w=0$, then 
$\frac{\partial}{\partial w}[\phi\circ(\phi^{(0)})^{-1}]=\tilde\gamma$.
Since $\eta_{\bar k}=(\det\phi)^{-1}(\tilde\eta_{\bar k})$, 
following the proof of Lemma~3.7 of \cite{MW}, we get
$$
\frac{\partial}{\partial w}(\eta_{\bar0}\otimes\eta_{\bar1}^{-1})=
-\sum_{k=0,1}(-1)^k\Tr(\tilde\gamma\tilde Q_{\bar k})
\,\eta_{\bar0}\otimes\eta_{\bar1}^{-1}.
$$
Therefore $\tau_\TT(X,H,r)$ is invariant under this deformation.
\end{proof}

\subsection{Behavior under $T$-duality}
Recall that if $\pi\colon X\to M$ is a principal $\TT$-bundle and $H$ is a $\TT$-invariant
closed $3$-form with integral periods, the $T$-dual bundle $\hat\pi\colon\hat X\to M$
is a principal bundle whose structure group is the dual circle group $\hat\TT$.
It is topologically determined by $c_1(\hat X)=\pi_*[H]$.
We now explain $T$-duality at the level of differential forms.
We have the commutative diagram
$$
\xymatrix@=2.5pc@ur{X\ar[d]_{\pi} & X\times_M\hat X\ar[d]^{\hat p}\ar[l]_p\\
M & \hat X\ar[l]^{\hat\pi}}
$$
where $X\times_M\hat X$ denotes the correspondence space.
Then $p^*[H]=\hat p^*[\hat H]\in H^3(X\times_M\hat X,\ZA)$.
Choosing connection 1-forms $A$ and $\hat A$ on the circle bundles $X$ and
$\hat X$, respectively, the formula
$$
T(\omega)=\int_\TT e^{p^*\!A\wedge\hat p^*\!\hat A}\,p^*\omega,\quad
\omega\in\Om^\bullet(X)
$$
gives linear map
$T\colon\Om^{\bar k}(X)\to\Om^{\overline{k+1}}(\hat X)$, $k=0,1$.
Similarly, we define $S\colon\Om^{\bar k}(X)\to\Om^{\overline{k+1}}(\hat X)$ by
$$
S(\hat\omega)=\int_{\hat\TT}e^{-p^*\!A\wedge\hat p^*\!\hat A}\,\hat p^*\hat\omega,
\quad\hat\omega\in\Om^\bullet(\hat X).
$$

We next explain the construction of the $T$-dual flux form $\hat H$ on $\hat X$.
Let $\pi^*F=dA$ and $\hat\pi^*\hat F=d\hat A$ be the curvatures of the 
connections $A$ and $\hat A$, respectively.
Since $H-A\wedge\hat F$ is a basic differential form on $X$, we have
$$
H=A\wedge\pi^*\hat F-\pi^*\Om
$$
for some $\Om\in\Om^3(M)$.
Define the $T$-dual flux $\hat H$ by
$$
\hat H=\hat\pi^*F\wedge\hat A-\hat\pi^*\Om.
$$
Then $\hat H$ is closed.
Since
$$
d(p^*\!A\wedge\hat p^*\!\hat A)=-p^*H+\hat p^*\hat H,
$$
we have
$$
T\circ d_H=d_{\hat H}\circ T,\quad d_H\circ S=S\circ d_{\hat H}.
$$
Therefore $T$-duality induces isomorphisms on twisted cohomology groups
$$
T_*\colon H^{\bar k}(X, H)\to H^{\overline{k+1}}(\hat X,\hat H),\quad k=0,1
$$
with inverse $S_*$ \cite{BEM}, and there is an isomorphism
$$
\det T_*\colon\det H^\bullet(X,H)\cong
(\det H^\bullet(\widehat X,\widehat H))^{-1}.
$$
We will relate the twisted analytic torsions under this identification.
 
Given the Riemannian metric $g_X$ on $X$ or a triple $(g_M,A,r)$, we
define the $T$-dual metric on $\hat X$ as
$$  
g_{\hat X}=\hat\pi^*g_M+r^{-2}\hat A\odot\hat A 
$$
or given by the triple $(g_M,\hat A,r^{-1})$ so that $g_{\hat X}$ is
$\hat\TT$-invariant and the length of every fiber is $r^{-1}$.
We study the $T$-duality map on invariant differential forms.

\begin{lemma}\label{Riem_X}
Under the above choices of Riemannian metrics, 
$$
T\colon\Omega^{\bar k}(X)^\TT\to\Omega^{\overline{k+1}}(\hat X)^{\hat\TT},
\quad
S\colon\Omega^{\bar k}(\hat X)^{\hat\TT}\to\Omega^{\overline{k+1}}(X)^\TT
$$
are isometries for $k=0,1$.
\end{lemma}

\begin{proof}
For any $\omega=r^{1/2}\pi^*\omega_1+r^{-1/2}A\wedge\pi^*\omega_2
\in\Om^{\bar k}(X)^\TT$,
$T(\omega)=r^{-1/2}\hat\pi^*\omega_2+r^{1/2}\hat A\wedge\hat\pi^*\omega_1$.
The result follows from applying formula in the proof of Lemma~\ref{isom} to
both $\omega$ and $T(\omega)$.
$S$ is the inverse of $T$.
\end{proof}

\begin{theorem}[$T$-duality and analytic torsion for circle bundles]
In the above notations, we have, up to a sign,
$$
(\det T_*)(\tau_\TT(X,H,r))=\tau_{\widehat \TT}(\widehat X,\widehat H)^{-1}
\in(\det H^\bullet(\widehat X,\widehat H,r^{-1}))^{-1}.
$$
\end{theorem}

\begin{proof}
We denote the restriction of $d_{\hat H}$ to $\Om^{\bar k}(\hat X)^{\hat\TT}$
by $\hat d_{\bar k}$.
Since $T$ is an isometry, we have
$T\circ d_{\bar k}^\dagger=\hat d_{\overline{k+1}}^\dagger\circ T$
and hence $T\circ(d_{\bar k}^\dagger d_{\bar k})=
(\hat d_{\overline{k+1}}^\dagger\hat d_{\overline{k+1}})\circ T$.
It follows that $T$ isometrically maps the space of $H$-twisted even (odd)
degree harmonic forms on $X$ to the space of $\hat H$-twisted odd (even)
degree harmonic forms on $\hat X$.
So $T$ maps the unit volume elements of $H^\bullet(X,H)$ to those of
$H^\bullet(\hat X,\hat H)$ up to a sign.
$T$ also maps isometrically on other eigenspaces, preserving the (positive)
eigenvalues. 
We deduce that
$$
\zeta_\TT(s,d_{\bar k}^\dagger d_{\bar k})=\zeta_{\hat\TT}
(s,\hat d_{\overline{k+1}}^\dagger\hat d_{\overline{k+1}})
$$
and the result follows.
\end{proof} 

\medskip\noindent
{\bf Acknowledgments.} We thank the referees for helpful comments.\\

\end{document}